\documentclass[11pt,leqno]{amsart}
\usepackage{graphicx}
\usepackage{geometry}
\usepackage{amsmath}
\usepackage{amssymb}
\usepackage{amsthm}
\usepackage{amsfonts}
\usepackage{enumerate}
\usepackage[scriptsize,hang,raggedright]{subfigure}
\usepackage[cmyk]{xcolor}
\usepackage{esint}

\usepackage{hyperref}
\usepackage{dsfont}
\usepackage{color}

\usepackage{cite}
\usepackage{bbm}
\usepackage{multirow}

\numberwithin{equation}{section} 

\newtheorem{theorem}{Theorem}[section]
\newtheorem{proposition}[theorem]{Proposition}
\newtheorem{lemma}[theorem]{Lemma}

\theoremstyle{remark}
\newtheorem{remark}[theorem]{Remark}
\newcommand{\R}{\mathbb{R}}

\newcommand{\ba}{\begin{array}}
\newcommand{\ea}{\end{array}}

\newcommand{\bthm}{\begin{theorem}}
\newcommand{\ethm}{\end{theorem}}
\newcommand{\bprop}{\begin{proposition}}
\newcommand{\eprop}{\end{proposition}}
\newcommand{\blemma}{\begin{lemma}}
\newcommand{\elemma}{\end{lemma}}

\newcommand{\beqn}{\begin{equation}}
\newcommand{\eeqn}{\end{equation}}
\newcommand{\beqns}{\begin{equation*}}
\newcommand{\eeqns}{\end{equation*}}

\newcommand{\supp}{\operatorname{supp}}

\newcommand{\pr}{\prime}

\newcommand{\arrow}{\rightarrow}

\newcommand{\rtwo}{\mathbb{R}^2}

\newcommand{\Rd}{\mathbb{R}^d}

\newcommand{\Prob}{\mathcal{P}}
\newcommand{\A}{\mathcal{A}}
\newcommand{\Scal}{\mathcal{S}}

\renewcommand{\leq}{\leqslant}
\renewcommand{\geq}{\geqslant}

\definecolor{mygreen}{rgb}{0.1,0.75,0.2}

\newcommand{\eps}{\epsilon}

\newcommand{\E}{\mathsf{E}}

\newcommand{\Om}{\Omega}

\title[Interaction of Sets]{Nonlocal Shape Optimization Via Interactions of Attractive and Repulsive Potentials}

\author{Almut Burchard}
\address{Department of Mathematics, University of Toronto, Toronto, ON Canada}
\email{almut@math.toronto.edu}
\author{Rustum Choksi}
\address{Department of Mathematics and Statistics, McGill University, Montr\'{e}al, QC Canada}
\email{rustum.choksi@mcgill.ca}
\author{Ihsan Topaloglu}
\address{Department of Mathematics and Statistics, McMaster University, Hamilton, ON Canada}
\email{itopalog@math.mcmaster.ca}

\date{\today}                                        

\begin{document}

\begin{abstract}
We consider a class of nonlocal shape optimization problems for sets
of fixed mass where the energy functional is given by
an attractive/repulsive interaction potential in power-law form. 
We find that the existence of minimizers of this shape optimization problem 
depends crucially on the value of the mass. Our results include 
existence theorems for large mass and nonexistence theorems for 
small mass in the class where the attractive part of the potential
is quadratic. In particular, for the case where the repulsion is given by the Newtonian potential,
we prove that there is a critical value for the mass, above
which balls are the unique minimizers, and below which minimizers
fail to exist. The proofs rely on a relaxation of
the variational problem to bounded densities, and recent 
progress on nonlocal obstacle problems.
\end{abstract}
\maketitle

%%%%%%%%%%%%%%%%%%%%%%%%%%%%%%%%%%%%%%%%%%%%%%%%%%%%%%%%%%%%%%%%%%%%%%%%%%%%%%%%%%%%%%%%%%%%%%%%%%%
%%%%% INTRODUCTION
%%%%%%%%%%%%%%%%%%%%%%%%%%%%%%%%%%%%%%%%%%%%%%%%%%%%%%%%%%%%%%%%%%%%%%%%%%%%%%%%%%%%%%%%%%%%%%%%%%%
\section{Introduction}
\label{sec:intro}

In this note we address the following 
 nonlocal shape optimization problem: 
\beqn\label{e:gp}\tag{\bf P}
		\begin{gathered}
		\text{Minimize}\qquad \E(\Om) \ 
:= \ \int_{\Om} \int_{\Om} K(x-y)\,dx\, dy \\
\text{over measurable sets }\Om\subset\Rd\ (d\geq2)\ 
\text{of finite measure }|\Om |=m.
		\end{gathered}
	\eeqn
Here $K:\Rd\to\R\cup\{+\infty\}$ is a locally integrable, lower semicontinuous, radial function, and $|\Om |$ denotes the Lebesgue measure of the set $\Om$. In particular, we are interested in interaction potentials in the power-law form
	\beqn
			K(x)\, :=\, \frac{|x|^{q}}{q}-\frac{|x|^{p}}{p}
		\label{e:kernel}
	\eeqn
where $-d<p<q$ with $p, q \neq 0$. 

\begin{figure}[ht!]
     \begin{center}
 \subfigure[{\footnotesize $-d<p<0$ and $q>1$}]{
            \label{fig:first}
            \includegraphics[width=0.30\linewidth]{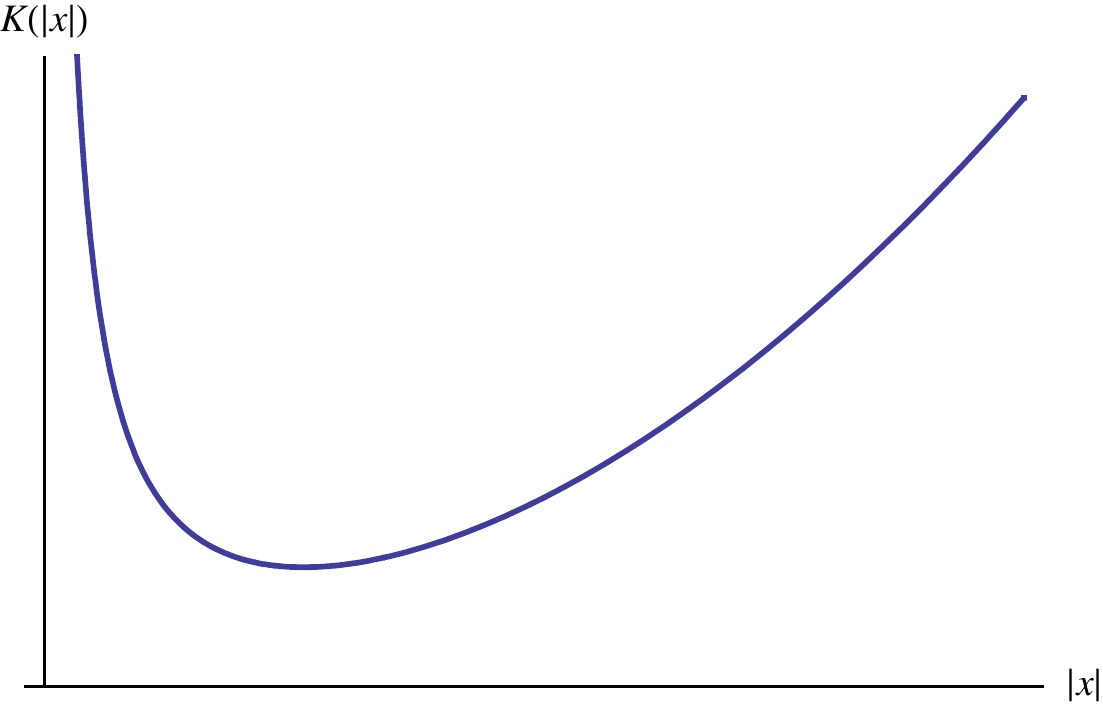}
        }\qquad\quad
        \subfigure[{\footnotesize $-d<p<0$ and $0<q<1$}]{
           \label{fig:second}
           \includegraphics[width=0.30\linewidth]{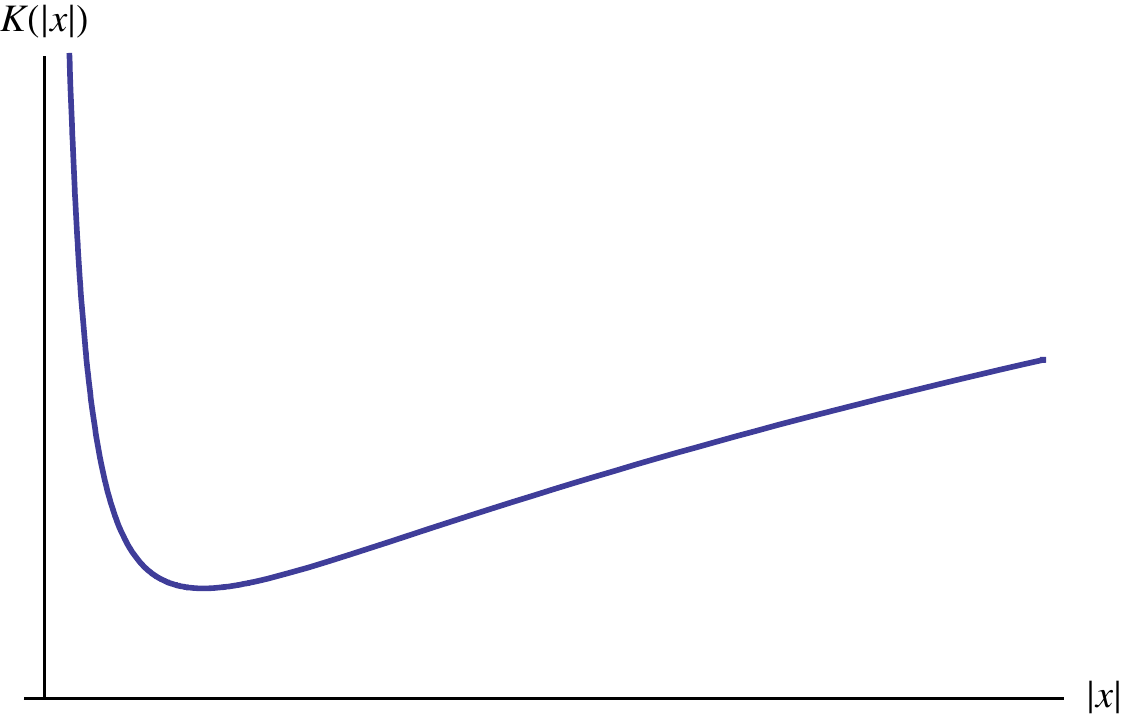}
        }
    \end{center}
    \caption{
        Generic examples of $K$ for various values of $q$ and $p$.
     }
   \label{Kdependingq}
\end{figure}

These sums of  
attractive and  repulsive power-law potentials have collective effect which is repulsive at short ranges but attractive at long ranges (see Figure \ref{Kdependingq}). 
We will focus on positive attraction $q>0$ and Riesz 
potential repulsions $-d < p <0$; 
the majority of our results pertain to quadratic attraction $q=2$, and some require $p$ to be at or below $2-d$, the exponent of the Newtonian potential. Our results are valid in any dimension $d \geq 2$ with the understanding that when $d=2$ the Newtonian repulsion (corresponding to $p=2-d=0$) is given by $\log|x|$, i.e., the kernel \eqref{e:kernel} is $K(x)=(1/q)|x|^q - \log(|x|)$ when $p=2-d$. Moreover, when $p=2-d$ the repulsive part of the energy is determined by the $H^{-1}$-norm of the characteristic function and is equal to $\|\chi_\Om\|_{H^{-1}}^2$ up to a constant. We use the notation $\chi_\Om$ for the characteristic (indicator) function of a set.

The problem \eqref{e:gp} is a toy example of
shape optimization problems where repulsive interactions at short distances
compete with attraction at long distances. As far as we know this is the first work to address such problems. 
It is closely related to the problem of minimizing
the nonlocal interaction energy
	\beqn\label{eqn:relax_energy}
		\E(\rho) :=
\int_{\Rd}\!\int_{\Rd} K(x-y)\,\rho(x)\rho(y)\,dx\,dy
	\eeqn
over non-negative densities $\rho\in L^1(\Rd)$ of
given mass $\|\rho \|_{L^1(\Rd)} \, = \, m$. Such 
functionals appear in a class of  well-studied 
self-assembly/aggregation models (e.g. see 
\cite{BernoffTopaz, FeHuKo11,HoPu2005,HuBe2010} 
%HoPu2006, LeToBe2009
and the references therein). Under broad assumptions
on the kernels, the existence of global 
minimizers \cite{CCP,CaChHu,ChFeTo14,SiSlTo2014} and qualitative properties 
of local minimizers \cite{Balague_etal13,CDM14} 
of these energies along with convex approximations of minimizers via analytical \cite{CrTo15} 
and numerical \cite{BaSh15} techniques have recently been investigated. These results
do not directly extend to \eqref{e:gp}, because 
a sequence of densities given by the indicator functions 
of measurable sets may converge weakly to densities 
taking values strictly between zero and one.
Nevertheless, we are able to exploit the relation between
the two problems to obtain existence and 
non-existence results for \eqref{e:gp}.

The purpose of this study is to lay out the foundations for addressing \eqref{e:gp},  
focusing mostly on the case of quadratic attraction.  In particular, we prove: 

\bthm
\label{thm:Coulomb} 
Let $K$ be of the form \eqref{e:kernel}, and let $m>0$.
\begin{itemize}
\item[(i)] 
For $q=2$ and $-d<p\leq 2-d$ and for sufficiently small $m$, the problem \eqref{e:gp} does not have a solution.
\item[(ii)]   For $q=2$ and $-d<p<0$ and for sufficiently large $m$, the ball of volume $m$ is the unique solution of \eqref{e:gp} up to translations.
\item[(iii)]  For $q=2$ and $p=2-d$, the unique solution  of \eqref{e:gp} is a ball of volume $m$ if $m\geq \omega_d$,  where $\omega_d$ denotes the volume of the unit ball in $\Rd$. If $m < \omega_d$,  the problem \eqref{e:gp} does not have a solution.
%\item[(iv)] For $q>0$ and $-d<p<2-d$ and for sufficiently small $m>0$, the problem \eqref{e:gp} does not have a solution.
\end{itemize}
\ethm

Our approach to Theorem \ref{thm:Coulomb} is via a {\it relaxation} of \eqref{e:gp} wherein the energy
\eqref{eqn:relax_energy} is minimized over
densities $\rho$ with $0\leq \rho\leq 1$ almost everywhere. We will denote this relaxed problem by \eqref{e:rp} and note that     
existence of minimizers was recently established in \cite{ChFeTo14}.
In Section \ref{sec:relax}, we  show  that \eqref{e:gp} has a solution if and 
only if the relaxed problem has a solution which is a characteristic function
(Theorem~\ref{thm:nec_suff_cond}). We also derive the first variation of \eqref{e:rp} and show that 
local minimizers are compactly supported.
These results hold for general kernels.
In Section \ref{sec:quad_attr} we turn our attention to power-law potentials and consider the quadratic attraction case. After establishing the uniqueness of minimizers we first prove part (i) of Theorem \ref{thm:Coulomb} via a recent regularity result of Carrillo, 
Delgadino and Mellet \cite{CDM14} for local minimizers of $\E$ over probability measures where they prove the connection with solutions of certain nonlocal obstacle problems and utilize their regularity \cite{Caf98,Sil2007}. Then we show that balls satisfy the first-order variational inequalities corresponding to \eqref{e:rp} when the mass is sufficiently large and prove parts (ii) and (iii) of Theorem \ref{thm:Coulomb}. Our results exploit the special nature (convexity) of the energy $\E$ for $q = 2$. We believe the basic approach to their proof should extend to all $q>0$. We address the challenges of such extensions in Section \ref{sect-otherq} and also mention when we can expect minimizers that are not necessarily balls.

Our conclusions and the consideration of \eqref{e:gp} are motivated by a number of old and new
shape optimization problems which we now describe in the physically most relevant case of three
dimensions.

%%%%%%%%%%%%%%%%%%%%%%%%%%%%%%%%%%%%%%%%%%%%%%%%%%%%%%%%%%%%%%%%%%%%%%%%%%%%%%%%%%%%%%%%%%%%%%%%%%%
%%%%% RELATED PROBLEMS
%%%%%%%%%%%%%%%%%%%%%%%%%%%%%%%%%%%%%%%%%%%%%%%%%%%%%%%%%%%%%%%%%%%%%%%%%%%%%%%%%%%%%%%%%%%%%%%%%%%
\section{Related Shape Optimization Problems}
\label{sec:problems}

We start with a problem of Poincar\'{e} on the shape 
of a fluid \cite{Poincare1}. Assuming 
vanishing total angular momentum, the total potential energy 
in a fluid body, represented by a set $\Omega \subset {\mathbb{R}}^3$,
	is given by
	\[
	    - \int_\Omega  \int_\Omega  \frac{C}{|x - y|} \, dx \, dy,
	\]
where $-C |x - y|^{-1}$ is the Newtonian potential 
resulting from the gravitational 
attraction between two points $x$ and $y$ in the fluid,
and $C>0$ is a physical constant.
After rescaling,
\emph{Poincar\'{e}'s variational problem} is given by
\[    \begin{cases}
		\text{Minimize}\qquad - 
\int_{\Omega}\!\int_{\Omega} \frac{1}{|x-y|}\,  dx  \, dy \\
\text{over measurable sets  }\Om\subset\R^3\text{ with  }|\Om|=m.\end{cases}
	\]
Poincar\'{e} asserted that among all shapes with fixed mass, 
the unique shape of lowest energy is a ball, and
proved this statement for sufficiently smooth sets.  He referred to previous work of Lyapunov but was critical of its incompleteness. 
It was not until almost a century later that all the details 
were sorted out by Lieb \cite{Li2} wherein the heart of the 
matter lies in the rearrangement ideas of Steiner for the 
isoperimetric inequality. These ideas are captured in the Riesz 
rearrangement inequality and its development (cf. \cite{LiLo, Bu}).
On the other hand, the maximum energy is not attained,
as by breaking up the shape and spreading out one can drive the 
energy to $0$. 
 
Another classical  variational problem with 
similar conclusions is the \emph{isoperimetric problem}: 
\[  
\begin{cases}
		\text{Minimize}\qquad {\rm perimeter}\,(  \Omega) \\
		\text{over sets  }\Om\subset\R^3\text{ of finite perimeter with  }|\Om|=m.
		\end{cases}
	\]
It is of course well-known that the only minimizers
are balls. Again, the maximum does not exist. 

The energies in both these problems are 
purely attractive in that they share an, albeit different, incentive 
for set elements to stay together. When these
are placed in direct opposition by subtracting the
energies, one obtains the \emph{nonlocal isoperimetric problem},
 which stated in 
dimension $d = 3$ is  
\beqn \label{e:nlip} \tag{\bf NLIP}
\begin{cases}
		\text{Minimize}\qquad {\rm perimeter} \, ({\Omega}) \ + \  \int_{\Omega}\!\int_{\Omega} \frac{1}{|x-y|} \, dx  \, dy \\
\text{over sets  }\Om\subset\R^3\text{ of finite perimeter with  }|\Om|=m.
\end{cases}
\eeqn
Here, the Newton potential $|x-y|^{-1}$
represents the electrostatic repulsion between two points $x$ and $y$,
and the double integral represents
the Coulomb energy of a uniform charge distribution
on $\Omega$.
The two terms are now in direct competition: 
Balls are \emph{best} (minimizers) for the first 
term but \emph{worst} (maximizers) for the second.
This functional was first introduced in \cite{CP1, CP2} in studying the small volume fraction asymptotics of the 
Ohta-Kawasaki functional. It was conjectured that there 
exists a critical mass $m_c$ such that minimizers are balls 
for $m \leq m_c$ and  fail to exist otherwise. There has 
recently been much work on the \eqref{e:nlip}
(see e.g. \cite{BoCr14,J,FL,KnMu2013,KnMu2014,LO1,LO2}).
To date what is known is that there exist two constants 
$m_1 \leq m_2$ such that 
	\begin{itemize}
		\item[(i)] balls are the unique minimizers if 
$m \leq m_1$, and 
		\item[(ii)] minimizers fail to exist if $m > m_2$. 
	\end{itemize}
It remains open whether or not $m_1 = m_2$. 
Thus the heuristic picture emerges that the perimeter completely 
dominates up to a critical mass, beyond which the Coulomb 
repulsion is strong enough to break sets apart. 

In the \eqref{e:nlip} the attraction, that is the incentive for the set to remain together, is via perimeter, a local quantity involving derivatives,
while the repulsion results from a pairwise interaction potential.
As such %the competition involves 
the short and long-range interactions 
%which 
are inherently different\footnote{Recently there has 
also been a significant interest in nonlocal set interactions 
via nonlocal derivatives
(see e.g. \cite{AmDeMa2010,CaRoSa2010,FFMMM} 
and in particular \cite{Va13} for a review). Here the repulsion is of Riesz-type and the attraction is created by
the interaction of a set $\Om$ with
its complement $\Om^c$. Specifically, the nonlocal energy considered in these works is given by
\[
	\int_{\Rd}\!\int_{\Rd} \frac{(\chi_{\Om}(x)-\chi_{\Om}(y))^2}{|x-y|^{d+s}}\,dxdy
\]
for some $\Om\subset\Rd$ and $0\leq s \leq 1$. There has also been interest in nonlocal set interactions via cross 
interaction of two phases (cf. \cite{BoKnRo15,CiDeNoPo15}). 
}. It is thus natural to 
consider problems where both attraction and repulsion
are dictated by pairwise interaction potentials in power-law form, for example, minimizers of  
\beqn
\frac12
 \int_{\Omega}\!\int_{\Omega} |x-y|^2\, dx\,dy \ + \ 
\int_{\Omega}\!\int_{\Omega} 
\frac{1}{|x-y|} \,  dx\,dy 
\label{e:quadratic-Coulomb}
\eeqn
over sets  $\Omega\subset \mathbb{R}^3$ with $|\Omega| \, = \, m$.
This is the special case of \eqref{e:gp}
with $q=2$, $p=-1$ and 
$d=3$. It can be viewed as toy problem for  the total potential energy of {\it spring-like} media which at short distances experience Coulombic repulsion and at longer distances experience the usual Hookean attraction. 
As in the \eqref{e:nlip}, balls are best for the first term
but worst for the second.  However the role of the mass $m$ is reversed according to the different scaling of the attractive 
and repulsive terms in \eqref{e:quadratic-Coulomb} with repulsion dominating for small $m$ and 
attraction dominating for large $m$. 
While in the \eqref{e:nlip} the lack of existence of minimizers 
is due to mass escaping to infinity, here it is due to oscillations.
Moreover, unlike for the \eqref{e:nlip}, 
here we can explicitly identify the critical threshold below which minimizers fail to exist and above which the unique minimizer is a  ball.

In this short paper we make a first step at addressing existence vs. nonexistence
for the general problem \eqref{e:gp}, depending on the mass parameter $m$.
Here there is a surprising lack of general mathematical tools:   
For controlling the attractive part of the interaction potential,
there is nothing like the well-developed regularity 
theory for minimal surfaces, which greatly benefited
the analysis of both the local and nonlocal isoperimetric problems, and recently, the  analysis of variational problems with nonlocal derivatives.   
On the other hand, the Riesz rearrangement inequality  
which was the key to solving Poincar\'{e}'s problem, goes in the wrong direction.

Finally, we remark that we only consider locally integrable kernels although kernels that are not locally integrable and appear in crystallization problems are of great interest from the point of view of the calculus of variations.

%%%%%%%%%%%%%%%%%%%%%%%%%%%%%%%%%%%%%%%%%%%%%%%%%%%%%%%%%%%%%%%%%%%%%%%%%%%%%%%%%%%%%%%%%%%%%%%%%%%
%%%%% MASS SCALING
%%%%%%%%%%%%%%%%%%%%%%%%%%%%%%%%%%%%%%%%%%%%%%%%%%%%%%%%%%%%%%%%%%%%%%%%%%%%%%%%%%%%%%%%%%%%%%%%%%%
\section{Mass Scaling}
\label{sec:prelim}

Throughout we consider nonlocal interaction energies \eqref{eqn:relax_energy} over three different classes:
\begin{itemize}
\item $\Scal_m := $ Characteristic functions of measurable sets $\Om\subset\Rd$ with $|\Om|=m$; 
\item 		$	\A_{m,M} := \left\{\rho\in L^1(\Rd)\cap L^\infty(\Rd) \colon 
\|\rho\|_{L^1(\Rd)}=m \text{ and } 0\leq\rho(x)\leq M\text{ a.e.}\right\}$; 
\item			$\Prob(\Rd) := \text{ probability measures over }\Rd$. 
		\end{itemize}	
With an abuse of notation we denote the energy by $\E$ over each class; however, we emphasize the dependence on the admissible class using the notation $\E(\Om)$, $\E(\rho)$ and $\E(\mu)$, respectively, when needed. Note that minimization over $\Scal_m$ is precisely our shape optimization problem ({\bf P}). Clearly $\Scal_m\subset\A_{m,1}$ and $\A_{m,1}$ is the weak closure of $\Scal_m$ in the weak $L^1$-topology. 
%Furthermore, $|\Om|$ denotes the Lebesgue measure of a set $\Om$ and $\omega_d$ is the volume of the unit ball in $\Rd$.

Over $\Prob(\Rd)$ the minimal energy scales differently than on $\Scal_m$ or $\A_{m,M}$. When we consider the nonlocal energy \eqref{eqn:relax_energy} 
over density functions $\rho\in L^1(\Rd)$,  the shape of  
 minimizers is independent of the mass $m$: The problem is
homogeneous  in $\rho$, that is
	\[
		\E(c\rho)=c^2 \E(\rho)
	\]
for any $c>0$.

On the other hand, for
\eqref{e:gp} this is not the case since the attractive and 
repulsive parts of the interaction energy scale differently under 
a dilation. To see this let us split the energy into
its attractive and repulsive parts, $\E=\E_q-\E_p$, where
\[
\E_q(\Om)=\frac{1}{q}\int_{\Om}\!\int_{\Om} |x-y|^q\,dx\,dy  \quad {\rm and} \quad  \E_p(\Om)=\frac{1}{p}\int_{\Om}\!\int_{\Om} |x-y|^p\,dx\,dy. 
\]
Given a measurable set $\Om\subset\Rd$ of volume $m$,
and $t>0$, the dilated set
	\[
		t\Om:=\{x\in\Rd\colon t^{-1}x\in\Om\}
	\]
has mass equal to $t^dm$. The attractive
and repulsive parts of the energy satisfy
\[
\E_q(t\Om) \ =\ t^{2d+q}\,\E_q(\Om)  \quad {\rm and} \quad
\E_p(t\Om) = t^{2d+p}\,\E_p(\Om)\,.
\]
Choosing $t=m^{-1/d}$ and replacing $\Om$ with $t\Om$, 
we see that \eqref{e:gp} 
is equivalent to minimizing
	\beqn
\E(t\Om) = m^{2+\frac{q}{d}}\,\, \E_q(\Om)\,\, -\,\,  m^{2+\frac{p}{d}}\,\, \E_p(\Om) \qquad \hbox{\rm over sets of volume $|\Om|=1$.}
\label{e:scaling}
	\eeqn

Since $p<0<q$, we see from \eqref{e:scaling} that for
sets of large mass the energy is dominated 
by attraction, whereas for small mass
it is dominated by repulsion.  The separate effects of each term are characterized by the following   well-known application of the Riesz rearrangement inequality.

\bprop
\label{lem:E_r}
For every non-zero $r>-d$ and each $m>0$, balls are the unique minimizers of
the energy
	\[
		\E_r(\Om) = \frac{1}{r}\int_{\Om}\!\int_{\Om} |x-y|^r\,dx\,dy
	\]
among measurable sets $\Om\subset\Rd$ of measure $m$.
There is no maximum; the supremum takes
the value $+\infty$ for $r>0$, and $0$ for $-d<r<0$.
\eprop

\begin{proof} 
Given a set $\Om\subset\Rd$ of measure $m>0$, let $\Om^*$ be
the open ball of the same measure centered at the origin.
Since the kernel $K_r(x)=\frac1r |x|^r$ is radially increasing,
it follows from the classical Riesz rearrangement inequality \cite{LiLo, Bu} 
that
\[
\E_r(\Om^*) \ \leq \ \E_r(\Om).
\]
(Note that the sign of the factor $\frac1r$ compensates 
for the change of monotonicity when $r<0$.)
Since $K_r$ is strictly increasing,
equality holds only if $\Om$ agrees with $\Om^*$ up to a translation and a set
of measure zero~\cite{Li2}, that is, if $\Om$ itself is a ball.

For the second statement, construct maximizing
sequences of sets $\{\Om_n\}_{n\geq 1}$, where
each $\Om_n$ is union of $n$ balls of mass $m/n$
whose pairwise distance exceeds $n$.
\end{proof}

\bigskip

In light of \eqref{e:scaling}, if
the mass is large, the attractive interaction dominates 
and we expect that balls are global minimizers
for \eqref{e:gp}. If the mass is small, the repulsion dominates and
we expect that minimizers fail to exist:  Rather, a 
minimizing sequence converges weakly to a density function taking on values strictly between $0$ and $1$. We now make these statements precise.

%%%%%%%%%%%%%%%%%%%%%%%%%%%%%%%%%%%%%%%%%%%%%%%%%%%%%%%%%%%%%%%%%%%%%%%%%%%%%%%%%%%%%%%%%%%%%%%%%%%
%%%%% RELAXED PROBLEM
%%%%%%%%%%%%%%%%%%%%%%%%%%%%%%%%%%%%%%%%%%%%%%%%%%%%%%%%%%%%%%%%%%%%%%%%%%%%%%%%%%%%%%%%%%%%%%%%%%%
\section{The Relaxed Problem}
\label{sec:relax}

We consider the following relaxation
of \eqref{e:gp}:
\beqn\label{e:rp}\tag{\bf RP}
\text{Minimize}\qquad \E(\rho)=\int_{\Rd}\!\int_{\Rd} K(x-y)\,\rho(x)\rho(y)\,dx\,dy \qquad
\text{over} \quad \A_{m,1}. 
\eeqn
In this section we will work with radially symmetric kernels $K(\cdot)$ which are 
\beqn\label{e:KP}
 \text{locally integrable, nonnegative, lower semicontinuous, and satisfy }  \, \lim_{|x|\to\infty} K(x) = \infty. 
 \eeqn
Note that this class of kernels include power-law potentials of the form \eqref{e:kernel}.

The following existence result was first proved for power-law potentials in \cite{ChFeTo14}. To obtain the existence of minimizers for more general kernels we can use the arguments in \cite[Theorem 3.1]{SiSlTo2014} and obtain that a minimizing sequence is tight. Then combining this with the arguments in \cite[Theorem 2.1]{ChFeTo14} we can conclude that a minimizing sequence is compact, i.e., has a convergent subsequence in the class of admissible functions $\A_{m,1}$.

\bprop[{\bf Existence of solutions}]
\label{prop:exist_relax}
Under the assumptions of \eqref{e:KP}, 
 the problem \eqref{e:rp} admits a solution for each $m>0$.
\eprop

We say that a function $\rho$ is a 
\emph{local minimizer} of $\E$ in $\A_{m,1}$ 
(in the $L^1$-topology), 
if $\E(\rho)\leq \E(\rho+\phi)$ for all 
$\phi\in L^1(\Rd)$ with $\|\phi\|_{L^1}<\delta$ and
$\rho+\phi \in\A_{m,1}$. Local minimizers satisfy the 
following necessary  condition.

\blemma \label{lem:first_var}
Let $\rho$ be a local minimizer of the 
energy $\E$ in $\A_{m,1}$. Then there exists a constant $\lambda>0$ such that
(except for $x$ in a set of measure zero),
	\beqn
			 K * \rho (x)\quad \begin{cases}
			 							= \lambda \quad &\text{if}\ \ 0<\rho(x)<1,\\
										\geq \lambda \quad &\text{if}\ \ \rho(x)=0, \\
										\leq \lambda \quad &\text{if}\ \  \rho(x)=1.
			 						 \end{cases}
		\label{eqn:first_var}
	\eeqn
\elemma

\begin{proof}
We proceed as in \cite[Lemma 4.1.2]{Friedman} and \cite[Lemma 1.2]{CiDeNoPo15}.
Let $\rho\in\A_{m,1}$ be a local minimizer of $\E$. 
We need to construct perturbations that
are nonnegative on $S_0:=\{x:\rho(x)=0\}$, 
nonpositive on $S_1:=\{x:\rho(x)=1\}$, and preserve
mass.  Let $\phi$ and $\psi\in L^1(\Rd)$
be compactly supported, bounded, nonnegative functions 
with $\phi=0$ a.e. in $S_1$, $\psi=0$ a.e. in $S_0$, and
	\beqn
\int_{\Rd}\phi(x)\,dx=\int_{\Rd}\psi(x)\,dx =1\,.
		\label{e:test_func_cond}
	\eeqn
Fix $\eps>0$, and define
\begin{align*}
\phi_\eps(x) &:= 
\frac{1}{\|\phi\,\chi_{\{1-\rho>\eps\}}\|_{L^1(\Rd)}}\, 
\phi(x)\chi_{\{1-\rho(x)>\eps\}}(x), \\
\psi_\eps (x) & := \frac{1}{\|\psi\,\chi_{\{\rho>\eps\}}\|_{L^1(\Rd)}}
\,\psi(x)\chi_{\{\rho(x)>\eps\}}(x).
\end{align*}

By construction,
$\rho + t(\phi_\eps - \psi_\eps)$ lies in
$\A_{m,1}$ and the perturbation is small
for sufficiently small values of $t>0$.
Since $\rho$ is a minimizer, it follows
that
	\[
	0 \leq \lim_{t\to 0^+} 
\frac{\E(\rho+t(\phi_\eps-\psi_\eps))-
\E(\rho)}{t} = 2\int_{\Rd} K*\rho(x)\, (\phi_\eps-\psi_\eps)(x)\,dx.
\]

Clearly, $\phi_\eps\to \phi$ and $\psi_\eps\to \psi$ as
$\eps\to 0$. By dominated convergence, we can pass to the limit as 
$\eps\to 0$ and obtain
	\beqn\label{e:first_var_ineq}
		\int_{\Rd} K*\rho(x)\,(\phi-\psi)(x)\,dx \geq 0.
	\eeqn
By density, \eqref{e:first_var_ineq} holds for all
nonnegative functions $\phi,\psi$ in $L^1(\Rd)$ with
$\phi(x)=0$ on $S_1$, $\psi(x)=0$ on $S_0$, and 
$\|\phi\|_{L^1(\Rd)}=\|\psi\|_{L^1(\Rd)}=1$.
Minimizing and maximizing separately over
$\phi$ and $\psi$, we obtain a constant
$\lambda\in\mathbb{R}$ such that
\[
\inf 
\left\{ \int_{\Rd} K* \rho(x) \, \phi(x)\, dx \colon \|\phi\|_{L^1(\Rd)}=1,\ \phi\geq 0,\ \text{and}\ \phi=0\ \text{a.e.}\ 
\text{on}\ S_1\right\} \geq \lambda
\]
and
\[
\sup \left\{ \int_{\Rd} K* \rho(x) \, \psi(x) \, dx \colon \|\psi\|_{L^1(\Rd)}=1,\ \psi\geq 0,\ \text{and}\ \psi=0\ \text{a.e.}\ \text{on}\ S_0\right\}
\leq \lambda.
\]
In particular, $\lambda > 0$ since so are $K$, $\rho$ and $\psi$. We conclude that $K*\rho\geq \lambda$ a.e. on $\{x: \rho(x)<1\}$,
and $K*\rho\leq \lambda$ a.e. on $\{x: \rho(x)>0\}$, as claimed.
\end{proof}

\medskip

One consequence of Lemma \ref{lem:first_var} 
is that the minimizers of $\E$ over $\A_{m,1}$ 
are compactly supported. This fact was established 
in \cite{CCP} for minimizers of $\E$ over $\Prob(\Rd)$; a more direct approach was
used in~\cite[Proposition 1.11]{CiDeNoPo15}. In our situation, the argument is simple and we present it here for the convenience of the reader.

\blemma\label{lem:comp_supp}
Under the assumptions of \eqref{e:KP}, every local minimizer for
\eqref{e:rp} in $\A_{m,1}$ has compact support.
\elemma

\begin{proof} 
By Lemma \ref{lem:first_var},
there exists a constant $\lambda$ such that 
$K*\rho \leq\lambda$ almost everywhere on the support 
of $\rho$.  Changing $\rho$ on a set of measure zero,
if necessary, we may assume that $K*\rho(x)\leq \lambda$ for
{\em all} $x$ with $\rho(x)>0$.

Let $R>0$ be large enough such that
\[
C_R \, := \, \int_{|y|<R} \rho(y)\, dy>0.
\]
Since $K$ and $\rho$ are nonnegative,
we have for $x\in\Rd$ that
\begin{align*}
K*\rho(x) & \geq  \int_{|y|<R}K(x-y)\, \rho(y)\, dy\\
& \geq C_R \,  \inf\bigl\{K(z) \colon |z|>|x|-R\bigr\}.
\end{align*}
%where $\omega_d$ is the volume of the unit ball.
Therefore
\[
\lim_{|x|\to\infty} K*\rho(x) = \infty,
\]
and the sub-level set $\{x: K*\rho\leq \lambda\}$
is bounded.
Since the sub-level set contains the support of $\rho$, 
the claim follows.
\end{proof}

\bigskip A useful consequence of Lemma~\ref{lem:comp_supp} is
that $K*\rho$ is continuous (since $K$ is locally
integrable).  We can now reduce the geometric
variational problem to the relaxed problem.

\bthm[{\bf Necessary and sufficient conditions
for existence of \eqref{e:gp}}]\label{thm:nec_suff_cond}
Let $K$ be a radially symmetric kernel satisfying \eqref{e:KP}. Then 
the problem \eqref{e:gp} has a solution 
$\Om\subset \Rd$ if and only if its characteristic function
$\chi_\Om$ is a solution of \eqref{e:rp}.
\ethm

\begin{proof}
We will show that
\beqn\label{rp-gp}
\inf_{|\Om|=m} \E(\Om) = \inf_{\rho\in \A_{m,1}} \E(\rho)
\eeqn
and establish a  relationship between the solutions
of the two variational problems.
The inequality $\geq$ is trivial from the
definition of the two variational problems:  the 
characteristic function $\chi_\Om$ of any
set $\Om\subset\Rd$ of measure $m$ lies
in $\A_{m,1}$. Similarly if $\chi_\Om$ is
a global minimizer for $\E$, then clearly $\Om$ is  global minimizer for \eqref{e:gp}.

Conversely, suppose that the global minimum of $\E$ over
$\A_{m,1}$ is not achieved by a characteristic function,
and fix a global minimizer $\rho$. 
By Lemma~\ref{lem:comp_supp}, $\rho$ has compact support. 
Choose a sequence of measurable sets $\{\Om_n\}_{n\geq 1}$ whose characteristic functions
$\rho_n=\chi_{\Om_n}$ converge to $\rho$ weakly in $L^1(\Rd)$.
To be specific, take a dyadic decomposition of $\Rd$ into
cubes of side length $2^{-n}$, and let the intersection
of $\Om_n$ with a given cube $Q$ be the centered
closed subcube of volume $\int_Q\rho(x)\,dx$.
By construction, $|\Om_n|=m$, and $\rho_n\in\A_{m,1}$.
Since $\rho$ has compact support, the sets $\Om_n$ are contained
in a common compact set.

Clearly, $\rho_n\rightharpoonup\rho$ weakly in $L^1(\Rd)$.
It follows from the local integrability of $K$ that
\[
\lim_{n\to\infty} K*\rho_n(x) = K*\rho(x)
\] for every $x\in\Rd$, that is,
$K*\rho_n$ converges pointwise to $K*\rho$. By dominated 
convergence,
$K*\rho_n \to K*\rho$ strongly in $L^1(\Rd)$. 
Using once more that $\rho_n\rightharpoonup \rho$,
we conclude that
\[
\E(\Om_n)\, = \, \int_{\Om_n} K*\rho_n \ dx \ \to \ 
\int_{\Rd} (K*\rho)\, \rho \,  dx \, = \, \E(\rho).
\]
In particular,
\[
\inf_{|\Om|=m} \E(\Om) \leq \E(\rho) = \min_{\rho\in \A_{m,1}}\E(\rho)\,,
\]
and $\{\Om_n\}$ is a minimizing sequence for \eqref{e:gp}.
Since $\E(\Om)>\E(\rho)$
for every $\Om\subset\Rd$, no minimizer exists.
\end{proof}

%%%%%%%%%%%%%%%%%%%%%%%%%%%%%%%%%%%%%%%%%%%%%%%%%%%%%%%%%%%%%%%%%%%%%%%%%%%%%%%%%%%%%%%%%%%%%%%%%%%
%%%%% QUADRATIC ATTRACTION
%%%%%%%%%%%%%%%%%%%%%%%%%%%%%%%%%%%%%%%%%%%%%%%%%%%%%%%%%%%%%%%%%%%%%%%%%%%%%%%%%%%%%%%%%%%%%%%%%%%
\section{The Case of $q=2$}
\label{sec:quad_attr}

In this section we specialize to kernels of the form \eqref{e:kernel} where the attractive term is quadratic, i.e., $q=2$. The key observation here is that \eqref{e:rp}
can be rewritten as a convex minimization problem in the parameter regime $q=2$ and $-d < p < 0$, hence, allowing us to conclude the uniqueness of minimizers of the relaxed problem.

\blemma\label{lem:quad_uniq}
For $q=2$ and $-d<p<0$, the solution of
problem \eqref{e:rp} is unique up to translation, and is given by
a radial function.
\elemma

\begin{proof}
Since the energy $\E(\rho)$ is translation 
invariant, without loss of generality, we assume that $\int_{\Rd} x\rho(x)\,dx=0$. Then  
	\[
		\E_q(\rho)\, =\, \frac{1}{2}\int_{\Rd}\!\int_{\Rd}|x-y|^2\rho(x)\rho(y)\,dx\,dy \,= \, m\int_{\Rd}|x|^2\rho(x)\,dx,
	\]
and the attractive part of the energy is linear in $\rho$.

On the other hand, when $-d<p<0$, the 
repulsive part of the energy
	\[
	-\E_p(\rho)\, =\, 
-\frac{1}{p}\int_{\Rd}\!\int_{\Rd}|x-y|^p\rho(x)\rho(y)\,dx\,dy 
	\]
is strictly convex over $\A_{m,1}$ since the Fourier transform of the kernel $-K_p(x)=-\frac{1}{p}|x|^{p}$ is strictly positive when $-d<p<0$ \cite[Corollary 5.10]{LiLo}.

Therefore the energy is strictly convex among all functions in $\A_{m,1}$ with zero first moments, and the solution of \eqref{e:rp} is unique up to translations.

Radial symmetry of the solution follows from the uniqueness and, due to its isotropic nature, the rotational symmetry of the energy $\E(\rho)$ around the center of mass of any $\rho\in\A_{m,1}$.
\end{proof}

\begin{remark}\label{rem:log_uniqueness}
For $x\in\mathbb{R}^2$ we take
		\[
			K(x)=\frac{1}{2}|x|^2 - \log|x|
		\]
when $p=2-d$, and the repulsive part of the energy is given by
		\[
			-\E_p(\rho) = -\int_{\rtwo}\!\int_{\rtwo} \log\Big(|x-y|\Big)\rho(x)\rho(y)\,dxdy = C \|\rho\|_{H^{-1}}^2.
		\]
Hence, the repulsion term is strictly convex and we still have the uniqueness of minimizers in the case $p=2-d$ when $d=2$. 	
\end{remark}

\medskip
%%%%%%%%%%%%%%%%%%%%%%%%%%%%%%%%%%%%%%%%%%%%%%%%%%%%%%%%%%%%%%%%%%%%%%%%%%%%%%%%%%%%%%%%%%%%%%%%%%%
%%%%% NONEXISTENCE OF SET MINIMIZERS
%%%%%%%%%%%%%%%%%%%%%%%%%%%%%%%%%%%%%%%%%%%%%%%%%%%%%%%%%%%%%%%%%%%%%%%%%%%%%%%%%%%%%%%%%%%%%%%%%%%
\subsection{Nonexistence for \eqref{e:gp} for small mass}
\label{sec:nonexist}

To prove the nonexistence of minimizers in the small mass regime we specialize to kernels of the form \eqref{e:kernel} with $q=2$ and $-d < p \leq 2-d$.  
This range of Riesz potentials share some important properties via their correspondence to 
%One such property is that the energies defined via these potentials are reflection positive \cite{BuCh15,FrLi2010}, i.e.,
%	\[
%		|\E_p(\Om^+)| + |\E_p(\Om^-)| \, \geq \, 2|\E_p(\Om)|
%	\]
%where $\Om^+$ and $\Om^-$ are two possible symmetrizations of $\Om$ at a single hyperplane.
%Another important property is that
the obstacle problem for $(-\Delta)^s$ with $s\in(0,1]$ 
%corresponding to Riesz potentials $|x-y|^p$ with $-d<p\leq 2-d$ 
which enjoys rather strong regularity features \cite{Caf98,Sil2007}. This connection between the obstacle problem and nonlocal interaction energies over $\Prob(\Rd)$ was recently exploited by Carrillo, Delgadino and Mellet \cite{CDM14} to obtain regularity of local minimizers with respect to the $\infty$-Wasserstein metric $d_{\infty}$.\footnote{For $\mu$, $\nu\in\Prob(\Rd)$ the $\infty$-Wasserstein metric is defined as
	\[
	 d_\infty(\mu,\nu) \, := \, \inf_{\pi\in\Pi(\mu,\nu)} \sup_{(x,y)\in\supp\pi} |x-y|,
	\]
where $\Pi(\mu,\nu):=\{\pi\in\Prob(\Rd\times\Rd) \colon \pi(A\times\Rd)=\mu(A) \text{ and } \pi(\Rd\times A)=\nu(A) \text{ for all } A\subset\Rd \}$.
} Although, a priori local minimizers in the $d_\infty$-topology are not comparable with the local minimizers in the $L^1$-topology the regularity result is true for \emph{global} minimizers independent of the topology. Here we rephrase their results for interaction potentials in power-law form \eqref{e:kernel} (cf. \cite[Remark 3.1]{CDM14}). 

\begin{lemma}[Theorems 3.4 and 3.10 in \cite{CDM14}] \label{lem:reg}
Let $K$ be given by \eqref{e:kernel}. Let $\mu\in\Prob(\Rd)$ be a local minimizer of $\E$ over $\Prob(\Rd)$ in the topology induced by $d_\infty$.
	\begin{itemize}
	 \item[(i)]  If $q>0$ and $p=2-d$, then $\mu$ is absolutely continuous with respect to the Lebesgue measure and there exists a function $\phi \in L^{\infty}(\Rd)$ such that $d\mu(x)=\phi(x)\,dx$.
	 \item[(ii)] If $q>0$ and $p<2-d$, then $\mu$ is absolutely continuous with respect to the Lebesgue measure and there exists a function $\phi \in C^{\alpha}(\Rd)$ for all $\alpha<1$ such that $d\mu(x)=\phi(x)\,dx$.
	\end{itemize}
\end{lemma}

\medskip

\begin{remark}[$L^\infty$-control on global minimizers] \label{rem:Linfty_bd}
In the parameter regime $q>0$ and $-d< p \leq 2-d$ we can still control the $L^\infty$-bound of a \emph{global} minimizer. In fact, \cite[Theorem 1.4]{CCP} implies that any global minimizer $\mu\in\Prob(\Rd)$ of $\E$ over $\Prob(\Rd)$ is compactly supported. This, in light of Lemma \ref{lem:reg}(ii), yields that the density function $\phi$ is in $L^\infty(\Rd)$.
\end{remark}

Using these results we can relate the $L^\infty$-bound of minimizers to the mass constraint $m$ via scaling which in turn enables us to obtain nonexistence of minimizers of the set energy $\E(\Om)$ when the mass is sufficiently small.

\begin{proof}[{\bf Proof of Theorem \ref{thm:Coulomb}(i):}]
Let $\mu\in\Prob(\Rd)$ be a global minimizer of $\E$ over $\Prob(\Rd)$. Such a minimizer exists by \cite[Theorem 1.4]{CCP} or \cite[Theorem 3.1]{SiSlTo2014} in the parameter regime $q=2$, $-d<p\leq 2-d$. By Lemma \ref{lem:reg} and Remark \ref{rem:Linfty_bd}, $\mu$ is absolutely continuous with respect to the Lebesgue measure with bounded density, i.e., there exists a constant $C>0$ such that $\|\mu\|_{L^\infty}<C$ with an abuse of notation.

Consider $\rho_m:=m\,\mu$. For $m>0$ sufficiently small we have that $\rho_m\in\A_{m,1}$. Now we claim that $\rho_m$ minimizes $\E$ over $\A_{m,1}$. To see this let $\phi\in\A_{m,1}$ be an arbitrary function and note that $(1/m)\phi\in\Prob(\Rd)$. Using the fact that $\mu$ minimizes $\E$ over $\Prob(\Rd)$ and the scaling of the energy $\E$ we have that
	\[
		\E(\rho_m)\,  =\,  m^2\,\E(\mu) \leq m^2 \, \E\left(\frac{1}{m}\phi\right) \, =\,  \E(\phi).
	\]

On the other hand, by Lemma \ref{lem:quad_uniq} and Remark \ref{rem:log_uniqueness}, $\rho_m$ is the unique minimizer of $\E$ over $\A_{m,1}$ in any dimension $d \geq 2$. For $m$ sufficiently small we have $\|\rho_m\|_{L^\infty(\Rd)}=m\,\|\mu\|_{L^\infty(\Rd)}\leq m\,C <1$. Hence, when $m$ is small $\rho_m$ is not a characteristic function of a set. Since it is the unique solution to the problem \eqref{e:rp} by Theorem \ref{thm:nec_suff_cond} the energy $\E$ does not admit a minimizer over measurable sets of measure $m$. 
\end{proof}

\medskip
%%%%%%%%%%%%%%%%%%%%%%%%%%%%%%%%%%%%%%%%%%%%%%%%%%%%%%%%%%%%%%%%%%%%%%%%%%%%%%%%%%%%%%%%%%%%%%%%%%%
%%%%% EXISTENCE OF SET MINIMIZERS
%%%%%%%%%%%%%%%%%%%%%%%%%%%%%%%%%%%%%%%%%%%%%%%%%%%%%%%%%%%%%%%%%%%%%%%%%%%%%%%%%%%%%%%%%%%%%%%%%%%
\subsection{Existence for \eqref{e:gp} for large mass}
\label{sec:exist}

We first note that \emph{heuristically} Lemma \ref{lem:first_var} and Theorem \ref{thm:nec_suff_cond}  should imply 
existence for  $m \geq \omega_d$ in the case of Newtonian repulsion $p = 2-d$ and quadratic attraction $q =2$. 
To see this formally, assume that any local minimizer of \eqref{e:rp} is continuous on its support and let
	\[
 		\Omega \,= \, \{ x \in \Rd \colon  0 < \rho (x) < 1\}
 	\]
for a local minimizer $\rho$. Suppose, for a contradiction, that $|\Omega|>0$. Since we assume that $\rho$ is continuous on its support, $\Omega$ is an open set. Lemma \ref{lem:first_var} implies there exists a constant $\lambda$ such that 
\[ K\ast \rho (x) \, = \, \lambda \qquad {\rm on} \,\, \Omega. \]
Taking the Laplacian of both sides, we find for all $x \in \Omega$, 
\begin{eqnarray*}
\Delta  K\ast \rho (x) & = & \frac{1}{2} \Delta \left( |\cdot|^2 \ast \rho \right)  (x)  \,\, + \,\, 
\frac{1}{d-2}\Delta \left( \frac{1}{|\cdot|^{d-2}} \ast \rho \right)  (x) \\
& = & d \int_{\Rd} \rho (y) dy \,\, - \,\, d \omega_d \rho (x)  \,=\, 0, 
\end{eqnarray*}
or 
\[ \frac{m}{w_d} \, = \, \rho (x). \]
Hence if $m \geq \omega_d$ we obtain a contradiction  unless the set $\Omega$ is empty. This shows that for $m \geq \omega_d$, every local minimizer of \eqref{e:rp} must be a characteristic function. 
By Theorem \ref{thm:nec_suff_cond}, this establishes existence of \eqref{e:gp} for $m \geq \omega_d$ and characterizes the minimizer.
We will shortly prove this result rigorously and show that this lower bound is sharp. 

We now turn to the full range of Riesz potentials, i.e., to the regime $-d<p<0$. To prove the existence of set minimizers for the energy $\E$ when the mass $m$ is sufficiently large we will first prove that the characteristic function of a ball is indeed a critical point of the relaxed problem \eqref{e:rp}.

\blemma[Large balls satisfy the necessary condition of Lemma \ref{lem:first_var}]\label{lem:ball_crit}
Let any $q>1$ and $-d<p<0$. For sufficiently large mass $m$,
the characteristic function
of a ball of mass $m$ is a critical point
for the energy $\E$ on $\A_{m,1}$.
\elemma

\begin{proof} We split the kernel
into its attractive and repulsive parts by defining $K_q \, := \, (1/q)|x|^q$ and $K_p \, := \, (1/|p|)|x|^p$ so that $K=K_q+K_p$. Let $R$ be the radius of the ball of mass
$m$.
Since $K_q$ and $K_p$ are radial,
so are $K_q*\chi_{B_R}$ and $K_p*\chi_{B_R}$.

Since $K_q$ is radially increasing, so is $K_q*\chi_{B_R}$.
For $|x|\geq R/2$, we can estimate the radial derivative by
	\beqn\label{eqn:conv_derv_est1}
		\begin{aligned}
			\left(\nabla\left(K_q*\chi_{B_R}\right)(x)
\cdot\frac{x}{|x|}\right) 
&=\int_{|y|\leq R} \left|x-y\right|^{q-2}
\left(x-y\right)\cdot \frac{x}{|x|} \,dy \\
	 &\geq C_q\,R^{d+q-1},
		\end{aligned}
	\eeqn
where  the constant
\[
C_q \,= \,  \inf_{t \geq \frac12} \int_{|y|\leq 1} |te_1-y|^{q-2}(t-y_1)\, dy
\]
is positive since $q>1$ and $e_1$ denotes a unit vector in $\Rd$. 
%{\bf Need a bit more justification. Note that the integrand
%is odd under reflection at the line $y_1=t$,
%allowing to cancel the negative part.}

Similarly, $K_p*\chi_{B_R}$ is a decreasing function of $|x|$,
and we estimate for $|x|\geq R/2$,
	\beqn\label{eqn:conv_derv_est2}
	\left(\nabla\left(K_p*\chi_{B_R}\right)
(x)\cdot\frac{x}{|x|}\right)\, \geq\,  -C_p\,R^{d+p-1}
	\eeqn
for some constant $C_p>0$.

Let $R$ be sufficiently large so that $C_q\,R^q > C_p\,R^p$.
Such a number $R$ exists since $p<q$.
From \eqref{eqn:conv_derv_est1} and \eqref{eqn:conv_derv_est2} 
we get that $(K_q+K_p)*\chi_{B_R}(x)$ is increasing in $|x|$ 
for $|x|\geq R/2$. Therefore
	\[
		K*\chi_{B_R}(x)\,  \geq \, \lambda_R := K*\chi_{B_R}(x)
\Bigg|_{|x|=R}
	\]
for $|x|\geq R$. 
Furthermore,
\beqn
		K*\chi_{B_R}(x) \, < \, \lambda_R
\label{e:inner}
\eeqn for $R/2 \leq |x| < R$. 

We need to show that \eqref{e:inner}
extends to $|x|<R/2$.  We first note that since both $K_q * \chi_{B_R}$ and $K_p * \chi_{B_R}$ are radially symmetric we have that 
\beqn
	\begin{aligned}
\lambda_R\, &= \, \int_{|y|\leq R} \frac{|Re_1 - y|^q}{q} + \frac{|Re_1 - y|^p}{|p|}\,dy \\
				   &= \, R^{d+q} \int_{|y|\leq 1} \frac{|e_1-y|^q}{q}\,dy + R^{d+p} \int_{|y|\leq 1} \frac{|e_1-y|^p}{|p|}\,dy \\
				   &= \,\tilde{C}_q\,R^{d+q} + \tilde{C}_p\, R^{d+p} 
	\end{aligned}
\label{e:lambda}
\eeqn
where $\tilde{C}_q=K_q * \chi_{B_1}(x) \Big|_{|x|=1}>0$ and $\tilde{C}_p=K_p * \chi_{B_1}(x) \Big|_{|x|=1}>0$.

Using the fact that $K_q*\chi_{B_R}$ is increasing in $|x|$ and $K_p*\chi_{B_R}$ is decreasing in $|x|$, we estimate
	\beqn
		\begin{aligned}
	(K*\chi_{B_R})(x) &\leq \, (K_q*\chi_{B_R})(x)\Bigg|_{|x|=R/2} 
+ (K_p*\chi_{B_R})(0) \\ 
&= \, \tilde{\tilde{C}}_q\,R^{d+q}+\tilde{\tilde{C}}_p\,R^{d+p},
		\end{aligned}
	\notag
\eeqn
where
	\[
		\tilde{\tilde{C}}_q \,:=\, K_q * \chi_{B_1}(x) \Bigg|_{|x|=1/2}.
	\]
Hence, $\tilde{\tilde{C}}_q<\tilde{C}_q$ as $K_q*\chi_{B_R}$ is radially
increasing. Comparing this inequality with 
\eqref{e:lambda}, we see that
\eqref{e:inner} also holds for $|x|\leq R/2$, if $R$ is sufficiently
large.
\end{proof}

\medskip

\begin{proof}[{\bf Proof of Theorem \ref{thm:Coulomb}(ii):}]
By Lemmas \ref{lem:quad_uniq} and \ref{lem:ball_crit} the function $\chi_{B(0,R)}$ with $R=(m/\omega_d)^{1/d}$ 
is up to translation the unique solution of \eqref{e:rp} for $q=2$ and $-d<p<0$ provided $m$ is sufficiently large. By convexity,
it must be a global minimizer.
\end{proof}

Finally, as we noted in the introduction, in the case of Coulomb repulsion, i.e., when $p=2-d$, the thresholds of mass for existence/nonexistence appearing in Theorems \ref{thm:Coulomb} (i) and (ii) coincide and can be computed explicitly. This provides the complete picture regarding the minimization of $\E$ either over $\Scal_m$ or $\A_{m,1}$ in this special regime.

\begin{proof}[{\bf Proof of Theorem \ref{thm:Coulomb} (iii)}]
Consider the relaxed energy $\E$ over $\A_{m,1}$, and let $\rho_R:=\chi_{B(0,R)}$ with $R=(m/\omega_d)^{1/d}$ and $\rho_1:=(m/\omega_d)\chi_{B(0,1)}$. Note that both $\rho_R$ and $\rho_1$ are in $\A_{m,1}$.

Using the fact that $(d-2)^{-1}\int_{B(0,R)}|x-y|^{2-d}\,dy = d\omega_d \Phi(x)$ where $\Phi(x)$ solves the equation $-\Delta\Phi=\rho_R$ on $\Rd$ we can explicitly compute that
	\beqn
		K*\rho_R(x) \, =\,  \begin{cases} 
							\frac{m-\omega_d}{2}|x|^2 + \frac{d\omega_d R^2}{2(d-2)} + \frac{dmR^2}{2(d+2)} &\text{ if }|x|\leq R, \\
							\\
							\frac{m}{2}|x|^2 + \frac{\omega_d R^d}{(d-2)}|x|^{2-d} + \frac{dmR^2}{2(d+2)} &\text{ if }|x|>R.
					  \end{cases}
		\nonumber
	\eeqn
This shows via \eqref{eqn:first_var} that $\rho_R$ is a critical point of $\E$ over $\A_{m,1}$ if and only if $m\geq \omega_d$. Then by Lemma \ref{lem:quad_uniq} we get that $\rho_R$ is the unique minimizer of $\E(\rho)$ if and only if $m\geq \omega_d$. On the other hand, when $m < \omega_d$ a simple calculation shows that $\E(\rho_1) < \E(\rho_R)$. Moreover, by \cite[Theorem 2.4]{ChFeTo14}, $\rho_1$ is the unique global minimizer of $\E$ over $\A_{m,1}$ when $m < \omega_d$. Hence, the result follows by Theorem \ref{thm:nec_suff_cond}.
\end{proof}

\medskip

\begin{remark}[Failure of minimality of balls in 2-dimensions]
For more singular repulsive powers in 2-dimensions, we can determine the threshold below which the ball fails to be the global minimizer of $\E(\rho)$ by explicit calculations. When $d=2$, $q=2$ and $-2<p<0$, the energy of a ball of radius $R=(m/\pi)^{1/2}$ is given by
	\[
		\E(\chi_{B(0,R)}) \,= \, \frac{\pi}{2}R^6 + \frac{2\pi^2\Gamma(2+p)}{(-p)\Gamma\left(2+\frac{p}{2}\right)\Gamma\left(3+\frac{p}{2}\right)}R^{4+p},
	\]
where $\Gamma$ denotes the $\Gamma$-function. The computation of the attractive part of the energy is trivial; the computation of the repulsive part is given in \cite[Corollary 3.5]{KnMu2013}. On the other hand,
	\[
		\E(R^2\chi_{B(0,1)})\, =\, \left(\frac{\pi}{2}+ \frac{2\pi^2\Gamma(2+p)}{(-p)\Gamma\left(2+\frac{p}{2}\right)\Gamma\left(3+\frac{p}{2}\right)}\right)R^4.
	\]
Thus, choosing $R_c$ so that
	\[
		\frac{\pi}{2}R_c^2 + \left(\frac{2\pi^2\Gamma(2+p)}{(-p)\Gamma\left(2+\frac{p}{2}\right)\Gamma\left(3+\frac{p}{2}\right)}\right)R_c^p \,\, > \,\, \frac{\pi}{2}+ \frac{2\pi^2\Gamma(2+p)}{(-p)\Gamma\left(2+\frac{p}{2}\right)\Gamma\left(3+\frac{p}{2}\right)},
	\]
and noting that $R_c<1$ we see that for any $R\leq R_c$ we have that
	\[
		\E(R^2\chi_{B(0,1)}) \, \leq \, \E(\chi_{B(0,R)});
	\]
hence, $\chi_{B(0,R)}$ is not a global minimizer of $\E$ over $\A_{m,1}$.
\end{remark}

%%%%%%%%%%%%%%%%%%%%%%%%%%%%%%%%%%%%%%%%%%%%%%%%%%%%%%%%%%%%%%%%%%%%%%%%%%%%%%%%%%%%%%%%%%%%%%%%%%%
%%%%% POSITIVE ATTRACTION REGIME
%%%%%%%%%%%%%%%%%%%%%%%%%%%%%%%%%%%%%%%%%%%%%%%%%%%%%%%%%%%%%%%%%%%%%%%%%%%%%%%%%%%%%%%%%%%%%%%%%%%
\section{The Regime of $q>0$}
\label{sect-otherq} 

As we noted before the quadratic attraction case is special as the attractive part of the energy either over $\A_{m,1}$ or $\Prob(\Rd)$ is linear in its argument when we fix the center of mass of competitors to zero. This allows us to conclude the uniqueness of solutions to \eqref{e:rp}. The uniqueness of minimizers is key to the existence of solutions to \eqref{e:gp} as we utilize this to conclude that any stationary state to \eqref{e:rp} has to minimize the energy $\E$ over $\A_{m,1}$. When $q \neq 2$, on the other hand, even though Lemma \ref{lem:ball_crit} shows that the balls are stationary states in the parameter regime $q>1$, $-d<p<0$ when $m>0$ is large, due to the possible lack of uniqueness of minimizers, we cannot conclude the existence of solutions to \eqref{e:gp} for large measure. Nevertheless, we believe that the problem \eqref{e:gp} admits a solution for large values $m>0$ when $q>1$ as the energy is dominated by the attractive term which is minimized by balls of measure $m$.

The uniqueness of minimizers is also an important ingredient in establishing nonexistence of solutions to \eqref{e:gp}. Indeed, it is the uniqueness of solutions to \eqref{e:rp} which allows us to conclude that any solution of \eqref{e:rp} can be written as $m\mu$ for some $\mu$ that minimizes $\E$ over $\Prob(\Rd)$. Intuitively, for small $m>0$, the $L^\infty$-bound in the problem \eqref{e:rp} is not active, and the morphology of minimizers should be the same as of those over $\Prob(\Rd)$. When $m>0$ is large, on the other hand, the $L^\infty$-bound becomes active and adds addition repulsive effects to the problem penalizing accumulations.

When $q>0$ and $-d<p \leq 2-d$, nonexistence of solutions to \eqref{e:gp} as in Theorem \ref{thm:Coulomb}(ii) would also be true if the $L^\infty$-bound found in Lemma \ref{lem:reg} and Remark \ref{rem:Linfty_bd} was uniform for \emph{any} measure minimizer $\mu$. In that case, the proof of Theorem \ref{thm:Coulomb}(ii) would translate almost verbatim to the power regime $q>0$, $-d<p\leq 2-d$. A result in this direction is the following.

\bprop \label{prop:small_balls_not_min}
Let $K$ be of the form \eqref{e:kernel}. Then for $q>0$, $-d < p < 0$, and for $m>0$ sufficiently small the ball of measure $m$ is not a solution of \eqref{e:gp}.
\eprop

\begin{proof}
We will proceed by contradiction. If $B(0,r_n)$ with $\omega_d\,r_n^d=1/n$ were solutions of \eqref{e:gp} with $m=1/n$ for any $n\in\mathbb{N}$ then the weak limit of the sequence $\rho_n = n\chi_{B(0,r_n)} \in \Prob(\Rd)$ would also minimize the energy $\E$ over $\Prob(\Rd)$. This follows by noting that for fixed $\mu$ that globally minimizes $\E$ over $\Prob(\Rd)$ we have that for sufficiently large $n\in\mathbb{N}$
	\[
		\E(\mu) \leq \E(\rho_n) = n^2\,\E(\chi_{B(0,r_n)}) \leq n^2\,\E(n^{-1}\mu) = \E(\mu).
	\] 
The second inequality follows from \eqref{rp-gp}. 
 Thus $\lim_{n\arrow\infty} \E(\rho_n) = \inf_{\mu\in\Prob(\Rd)}\E(\mu)$, i.e., $\{\rho_n\}_{n\in\mathbb{N}}$ is a minimizing sequence for the energy $\E$ over $\Prob(\Rd)$. Arguing as in \cite[Theorem 3.1]{SiSlTo2014} via Lions' Concentration Compactness Theorem we obtain that $\rho_n$ has a weakly convergent subsequence and by the weak lower semicontinuity of $\E$ its limit minimizes $\E$ over $\Prob(\Rd)$. However, as $n\to\infty$, $\{\rho_n\}_{n\in\mathbb{N}}$ converges weakly to $\delta_0$, the Dirac measure at $x=0$,which has infinite energy. 
\end{proof}

A possible way of generalizing this result to conclude nonexistence  of  \eqref{e:gp}  for small $m$  is via the  \emph{energy-per-particle-pair}
	\beqn \label{eqn:energy_per_part}
		\eta(m) := \inf_{\rho\in\A_{m,1}} \frac{\E(\rho)}{m^2}
	\eeqn
associated with \eqref{e:rp}.  Because of the positivity of $K$, it is easy to see that if \eqref{e:gp} admits a solution for all $m>0$, then $\eta(m)$ is nondecreasing in $m$.  Moreover, if $\eta(m)$ is strictly increasing in $m$ (which is true when $q=2$, $-d<p<0$) then we would have the following  sufficient condition for nonexistence of minimizers: If $\eta^\pr(m_c)=0$ for some $m_c>0$, then  \eqref{e:gp} does not have a solution for $m<m_c$. 
Together with  Lemma \ref{lem:reg} and Remark \ref{rem:Linfty_bd}, this would prove nonexistence of  \eqref{e:gp} for sufficiently small $m>0$ when $q>0$ and $-d< p \leq 2-d$.  These remarks highlight  the fact that the  (strict) monotonicity of $\eta$ determines whether the $L^\infty$-constraint in $\A_{m,1}$ is active for the given value of $m$.

Finally, it remains to be proved whether there exists a regime of $m$, $q$ and $p$ where the minimizers are not balls. When $q$ is sufficiently large we expect that solutions to \eqref{e:gp} are rings rather than balls. Formally, the sequence of energies $\{\E(\rho)\}_{q>0}$ converges to
	\[
		\E_\infty (\rho) = \begin{cases}
									-\frac{1}{p} \int_{\Rd}\!\int_{\Rd} |x-y|^p\rho(x)\rho(y)\,dxdy &\mbox{if } \text{diam}(\supp\rho) \leq 1, \\
									+\infty				  &\mbox{otherwise}
								\end{cases}
	\]
as $q\to\infty$. Due to the purely repulsive effects in the energy $\E_\infty$ its minimizers $\rho$ should have convex supports and accumulate on the boundary of $\supp\rho$; however, these questions are open even in the Newtonian case.

\bigskip
%%%%%%%%%%%%%%%%%%%%%%%%%%%%%%%%%%%%%%%%%%%%%%%%%%%%%%%%%%%%%%%%%%%%%%%%%%%%%%%%%%%%%%%%%%%%%%%%%%%
%%%%% ACKNOWLEDGEMENTS
%%%%%%%%%%%%%%%%%%%%%%%%%%%%%%%%%%%%%%%%%%%%%%%%%%%%%%%%%%%%%%%%%%%%%%%%%%%%%%%%%%%%%%%%%%%%%%%%%%%
\noindent {\bf Acknowledgements.} The authors would like to thank the reviewer for their detailed reading and comments.
AB and RC were supported by NSERC (Canada) Discovery Grants. IT was supported by a Fields-Ontario Postdoctoral Fellowship.

%%%%%%%%%%%%%%%%%%%%%%%%%%%%%%%%%%%%%%%%%%%%%%%%%%%%%%%%%%%%%%%%%%%%%%%%%%%%%%%%%%%%%%%%%%%%%%%%%%%
%%%%% BIBLIOGRAPHY
%%%%%%%%%%%%%%%%%%%%%%%%%%%%%%%%%%%%%%%%%%%%%%%%%%%%%%%%%%%%%%%%%%%%%%%%%%%%%%%%%%%%%%%%%%%%%%%%%%%
\bibliographystyle{plain}
\bibliography{biblio}

\end{document}